\documentclass[11pt]{amsart}

\usepackage{amsmath,amssymb,amsthm}
\usepackage{thmtools,enumerate}
\usepackage{mathtools}

\usepackage[german,english]{babel}
\usepackage[autostyle]{csquotes}

\usepackage[all]{xy}
\usepackage{pstricks}

\usepackage{tikz-cd}
\usetikzlibrary{babel}

\usepackage{hyperref}
\hypersetup{
	colorlinks=true,
	linkcolor=red,
	citecolor=blue}

\theoremstyle{plain}

\newtheorem{thm}{Theorem}[section]
\newtheorem*{MainProb}{Main Problem}
\newtheorem{pro}[thm]{Proposition}
\newtheorem{lem}[thm]{Lemma}

\theoremstyle{definition}

\newtheorem{rem}[thm]{Remark}

\newcommand{\N}{\mathbb{N}}
\newcommand{\Q}{\mathbb{Q}}
\newcommand{\R}{\mathbb{R}}
\newcommand{\C}{\mathbb{C}}
\newcommand{\OO}{\mathcal{O}}
\newcommand{\mcal}{\mathcal}

\DeclareMathOperator{\NEb}{\overline{\mathrm{NE}}}
\DeclareMathOperator{\ddiv}{div}
\DeclareMathOperator{\Div}{Div}
\DeclareMathOperator{\Supp}{Supp}
\DeclareMathOperator{\mult}{mult}
\DeclareMathOperator{\Proj}{Proj}
\DeclareMathOperator{\Nef}{Nef}

\begin{document}
	
	\title[Programming the Minimal Model Program: a proposal]{Programming the Minimal Model Program:\\ a proposal}
	
	\author[V.\ Lazi\'c]{Vladimir Lazi\'c}
	\address{Fachrichtung Mathematik, Campus, Geb\"aude E2.4, Universit\"at des Saarlandes, 66123 Saarbr\"ucken, Germany}
	\email{lazic@math.uni-sb.de}

	\thanks{I gratefully acknowledge support by the Deutsche Forschungsgemeinschaft (DFG, German Research
Foundation) – Project-ID 286237555 – TRR 195. The idea for this paper originated during the Winter School on Computational Geometry at the Fraunhofer-ITWM in November-December 2022. I would like to thank J.\ Böhm and F.-O.\ Schreyer for extensive discussions on the implementation of the MMP over the years, to M.\ Musta{\c{t}}{\u{a}} and B.\ Mirgain for long discussions on Section \ref{sec:conedecomp}, and to G.\ Brown, W.\ Decker and the referees for very useful comments.
	\newline
		\indent 2020 \emph{Mathematics Subject Classification}: 14E30, 14Q15, 14Q20.\newline
		\indent \emph{Keywords}: Minimal Model Program, computer algebra, \texttt{OSCAR}}
	
	\begin{abstract}
	The aim of this paper is to propose a strategy to implement the Minimal Model Program in modern computer algebra systems.
	\end{abstract}

	\maketitle
	
	\begingroup
		\hypersetup{linkcolor=black}
		\setcounter{tocdepth}{1}
		\tableofcontents
	\endgroup
	
\section{Introduction}

The aim of this paper is to propose a strategy to programme the Minimal Model Program (MMP) in modern computer algebra systems such as \texttt{OSCAR} \cite{Oscar}.

The MMP is a programme in higher dimensional birational geometry aiming to classify varieties and pairs with mild singularities up to birational equivalence: the goal is to find in each such birational class at least one representative which is minimal in a very precise sense. When completed, the programme will show that -- up to birational equivalence -- all varieties with mild singularities are built out of three basic building blocks: varieties whose canonical class is either ample, numerically trivial or anti-ample. In other words, the theory predicts that one should be able to slice any such variety up into pieces whose curvature has almost everywhere the constant sign. Stated this way, the MMP bears resemblance to Thurston's Geometrisation in geometric topology.

The MMP has seen decisive progress ever since it was introduced by Mori in the 1980s. In particular, it is now complete for threefolds, and for varieties of log general type in every dimension. A comprehensive survey of the state of the art (as far as relevant for the present paper) can be found in \cite{Laz13}.

\medskip

In this paper I address the problem of possible \emph{implementations} of the Minimal Model Program on a computer. This is an important issue: we have very few concrete examples in higher dimensional geometry, see for instance the introduction to \cite{LS22}, and calculating minimal models of \emph{concrete} projective varieties is obviously a pressing issue, both from the point of view of applications within the MMP itself as well as applications in wider algebraic geometry.

There has been very limited progress on the problem of computer algebra implementations of the MMP. Ideally, at some point in the future one would feed a computer algebra system with a projective variety (given, for instance, by global equations; or by local equations in charts, requiring massive parallelisation \cite{BDFKPRR21}), tell it to calculate its minimal models, which then the system would spit out after a certain time. In order to achieve the implementation in this generality in a foreseeable future, a computer algebra system will need human input, simplifications and theoretical improvements at several steps along the way.

My goal in this paper is to propose a possible strategy for such an implementation based on the approach from \cite{CL12a,CL13,KKL16}; surveys of these ideas can be found in \cite{CL12b,Laz13}. The expectation is that the algorithm proposed here is realisable in modern or future computer algebra systems such as \texttt{OSCAR}. The input of the algorithm is the set of multidegrees of a set of generators of a certain finitely generated multigraded ring on a klt pair $(X,\Delta)$ of log general type, see Section \ref{sec:graded}; the output is a minimal model of $(X,\Delta)$ or, moreover, a sequence of maps in a $(K_X+\Delta)$-MMP which terminates with a minimal model of $(X,\Delta)$.

I will concentrate in this paper only on \emph{pairs of log general type} over $\C$, since this is the class of objects where the MMP is not a conjecture, but a theorem.

\section{Countable vs.\ uncountable problems}

Many problems in algebraic geometry are what I prefer to call \emph{uncountable} problems: problems whose decidability seems to depend on uncountably many parameters. A part of important theoretical progress on some of them is to translate them into \emph{countable} problems, i.e.\ problems whose decidability depends on at most countably many parameters. Then, at least in theory, a computer has a chance to verify a task in a finite time.

First I give two examples to illustrate the issue.

\subsection*{Example 1: Uniruledness}

A complex variety $X$ is \emph{uniruled} if it can be covered by (singular) rational curves; a good source on uniruledness is \cite{Deb01}, whereas a more comprehensive but more technical presentation is in \cite{Kol96}. Since this is a birational property, we may assume that $X$ is smooth. Deciding whether $X$ is uniruled seems a priori to be extremely difficult: one has to test whether through each point on $X$ there is a rational curve. One can relax this to testing whether through each \emph{general} point there is a rational curve, but practically this is not an improvement.

One of the main results of \cite{BDPP} is the proof that a smooth variety $X$ is uniruled if and only if the canonical class $K_X$ is not pseudoeffective, i.e.\ if the numerical class of $K_X$ is not a limit of the numerical classes of effective divisors. Assuming the full Minimal Model Program (including the Abundance conjecture) this is then equivalent to $K_X$ not being effective itself. In other words, we expect:

\begin{center}
$X$ uniruled$\quad\Longleftrightarrow\quad H^0(X,mK_X)=0$ for all $m>0$.
\end{center}

Note that this equivalence holds unconditionally in dimension $3$. Therefore, testing uniruledness is a countable problem and could in principle be implemented on a computer.

In practice this means that it suffices to find one positive integer $m$ such that $H^0(X,mK_X)\neq0$ in order to show that $X$ is not uniruled.

\subsection*{Example 2: Rational connectedness}

A complex variety $X$ is \emph{rationally connected} if through every two points on $X$ there is a (singular) rational curve; again, good sources on rational connectedness are \cite{Deb01,Kol96}. We may again assume that $X$ is smooth. A priori, testing whether $X$ is rationally connected seems to be an even more difficult problem than testing uniruledness. 

However, a conjecture attributed to Mumford suggests that testing rational connectedness is in fact a countable problem and could in principle be implemented on a computer. Indeed, one should have:

\begin{center}
$X$ rationally connected$\quad\Longleftrightarrow\quad H^0\big(X,(\Omega_X^1)^{\otimes m}\big)=0$ for all $m>0$.
\end{center}

This equivalence holds unconditionally in dimension $3$ by \cite{KMM92} and in almost all cases in dimension $4$ by \cite{LP17}.

In practice this means that it suffices to find one positive integer $m$ such that $H^0\big(X,(\Omega_X^1)^{\otimes m}\big)\neq0$ in order to show that $X$ is not rationally connected.

\subsection*{The MMP}

The Minimal Model Program has for a long time looked like an uncountable problem. The classical strategy for varieties $X$ with mild singularities goes like this:

\begin{enumerate}[\normalfont (i)]
\item first, by the Cone theorem \cite[Theorem 3.7]{KM98} there are at most countably many $K_X$-negative extremal rays in Mori's cone of curves $\NEb(X)$, which one can use to construct a contraction morphism to another variety $Y$. In the classical case of surfaces, one wants to find a $({-}1)$-curve and construct the corresponding morphism by Castelnuovo's contraction criterion. Finding such a curve seems to me to be an a priori uncountable problem and I do not see how it can be implemented on a computer.

\item Once one constructs such a morphism as in (i), if it is not birational, one stops. If it is birational, then one must decide whether it is a divisorial contraction (i.e.\ the good case) or a small contraction (i.e.\ the bad case). Equivalently, one must be able to test whether the canonical class on the target is $\Q$-Cartier or not, see \cite[2.6]{KM98}. If one is in the bad case, then one must additionally construct a flip \cite[Definition 2.8]{KM98}. The construction of the flip is potentially implementable on a computer.

\item One repeats (i) and (ii) until either one stops as in (ii), or the canonical class becomes nef. Testing whether a divisor (or a line bundle) on a variety is nef (i.e.\ whether it intersects every curve non-positively) seems to me to be an uncountable problem, similarly as in (i). Additionally, at no stage of the algorithm does one have control over how many repetitions one has to execute until the algorithm stops.
\end{enumerate}

\subsection*{The MMP with scaling}

One of the obvious issues in (i) above is that from the Cone theorem one does not know how to generate a suitable extremal ray and hence a contraction morphism. Instead, one can employ Shokurov's MMP \emph{with scaling} of some ample divisor \cite[Remark 3.10.10]{BCHM}: one fixes an ample (or just big) divisor $H$ on $X$ such that $K_X+H$ is nef, and determines the number
$$\lambda:=\inf\{t\in\R\mid K_X+tH\text{ is nef}\,\}.$$
Then one can show that $\lambda$ is in fact a rational number, that there exists a $K_X$-negative extremal ray $R\subseteq\NEb(X)$ such that $(K_X+\lambda H)\cdot R=0$, and that the divisor $K_X+\lambda H$ is semiample and defines a morphism which contracts precisely the curves whose classes belong to $R$. Then one continues with the steps (ii) and (iii) as above. This \emph{MMP with scaling of $H$} is much less arbitrary than a general MMP and it seems to be a step closer to being a countable problem: the main point is finding such $\lambda$, and for this one has to be able to calculate finitely many intersection products, and then repeat the procedure. However, the issue of testing whether a divisor is nef remains: to me this looks like an uncountable problem.

\subsection*{Another point of view}

A new outlook on the MMP was proposed in \cite{CL12a,CL13}: that once one knows that a certain graded ring is finitely generated, then the MMP becomes a problem closely related to the convex geometry of certain cones of divisors naturally associated to that graded ring. The goal of the remainder of the present paper is to show how this strategy can be adapted to an algorithm that can be implemented in a computer algebra system.

\section{Step 1: Graded rings and the Main Problem}\label{sec:graded}

In this paper, I consider the Minimal Model Program for projective $\Q$-factorial klt pairs $(X,\Delta)$ defined over $\C$ of \emph{log general type}, i.e.\ such that $K_X+\Delta$ is big. The canonical source on the foundational material in the MMP is \cite{KM98}. However, everything works in the category of varieties where analogous birational procedures can be executed as in \cite[Section 5]{KKL16}.

To start with, let $X$ be a projective $\Q$-factorial variety and let $D$ be a $\Q$-divisor on $X$. We define the global sections of $D$ by 
\[
H^0(X,D)=\{f\in k(X)\mid \ddiv f+D \geq 0 \}.
\]
Then for $\Q$-divisors $D$ and $D'$ on $X$ there is a well-defined multiplication map 
$$H^0(X,D)\otimes H^0(X,D')\to H^0(X,D+D'),$$
so that, if we are given a bunch of $\Q$-divisors $D_1,\dots,D_r$ on $X$, we can define the corresponding \emph{divisorial ring} as
\[
\mathfrak R=R(X;D_1, \dots, D_r)=\bigoplus_{(n_1,\dots, n_r)\in \N^r} H^0(X,
n_1D_1+\dots + n_rD_r),
\]
and the corresponding cone 
$$\mcal C=\sum_{i=1}^r\R_+ D_i\subseteq \Div_\R(X).$$
If $\mathfrak R$ is finitely generated, then we can define the \emph{support} of $\mathfrak R$, denoted by $\Supp\mathfrak{R}$: this is the convex hull of all integral divisors $D\in\mcal C$ such that $H^0(X,D)\neq0$. It is easily seen that $\Supp\mathfrak R$ is a rational polyhedral cone: indeed, pick finitely many generators $f_i$ of $\mathfrak R$, and let $E_i\in\mcal C$ be the divisors such that $f_i\in H^0(X,E_i)$. Then $\Supp\mathfrak R=\sum\R_+E_i$.

The following result gives the most important example of a finitely generated divisorial ring.

\begin{thm}\label{thmA}
Let $X$ be a normal projective variety and let $\Delta_i$ be $\Q$-divisors on $X$ such that each $\Delta_i$ is big and each pair $(X,\Delta_i)$ is klt for $i=1,\dots,r$. Then the ring
$$ R(X;K_X+\Delta_1,\dots,K_X+\Delta_r) $$
is finitely generated.
\end{thm}

This was first proved in \cite[Corollary 1.1.9]{BCHM} by the full machinery of the MMP, and then without the MMP in \cite[Theorem 1.2]{Laz09} and \cite[Theorem A]{CL12a}; see also \cite[Theorem 2]{CL13} for the above formulation. In particular, this result implies the finite generation of the \emph{canonical ring} $R(X,K_X)$ on a smooth projective variety $X$, which was implicitly conjectured in Zariski's famous paper \cite{Zar62}.

It is one of the main open problems in higher dimensional birational geometry to prove an analogue of Theorem \ref{thmA} when the divisors $K_X+\Delta_i$ are not big -- note that this is indeed a theorem in dimensions at most $3$. The case of threefolds of not necessarily general type is probably the most important case for the implementation of the MMP on a computer, beyond the case of pairs of log general type.  

\medskip

The first step in the potential algorithm of the implementation of the MMP is the following:

\begin{MainProb}\label{prob1}
Let $(X,\Delta)$ be a projective $\Q$-factorial klt pair of log general type, and let $A_1,\dots,A_r$ be ample divisors on $X$. Find the \textbf{multidegrees of some set of generators} of the ring
\[
R\big(X;m_0(K_X+\Delta),m_1(K_X+\Delta+A'_1),\dots,m_r(K_X+\Delta+A'_r)\big)
\]
for some positive integers $m_i$, and some effective divisors $A_i'\equiv A_i$.
\end{MainProb}

As we will see later, it is important for the MMP algorithm in this paper that the ample divisors $A_1,\dots,A_r$ are chosen in such a way that their numerical classes generate the N\'eron-Severi space $N^1(X)_\R$.

\medskip

The Main Problem is clearly a countable problem: indeed, one might be able to find a brute-force algorithm to calculate the graded pieces and test the generation relations in the ring -- since each graded piece is a finite dimensional vector space, one has at most countably many generators and relations. However, such a brute-force algorithm will not in general end in a finite time without an algebro-geometric input about the specific examples that one wants to calculate on a computer.

\medskip

\emph{However}: 

\subsection*{The Main Problem on surfaces}

Finding an \emph{explicit bound} of the degree of generation as in the Main Problem -- at least in explicit examples -- is one of the technical prerequisites for this algorithm to work.

There is evidence that the degree of generation of the divisorial ring as in the Main Problem is bounded uniformly, depending only on the dimension and the number of prime divisors in the supports of the divisors $K_X+\Delta$ and $A_i$. A version of this holds on surfaces by \cite[Proposition 4.1]{CZ14} and \cite[Theorem 3.1]{CL14} and is conjectured to hold in any dimension.

More precisely, we have \cite[Theorem 3.1]{CL14}:

\begin{thm}\label{t_efgs}
Let $X$ be a smooth surface and let $S_1,\dots,S_p$ be distinct prime divisors such that $\sum_{i=1}^p S_i$ is a simple normal crossings divisor. Let $B_0,B_1,\dots,B_r$ be $\Q$-divisors on $X$ such that:
\begin{enumerate}[\normalfont (i)]
\item $\Supp B_j=\sum_{i=1}^p S_i$ for all $j$, and
\item $\lfloor B_j\rfloor=0$ for every $j$, 
\end{enumerate}
and let $k$ be a positive integer such that all $kB_j$ are Cartier. Then there exists a positive integer $m$ depending only on $p $ and $k$ such that the ring 
$$R\big(X;m(K_X+B_0),m(K_X+B_1),\dots,m(K_X+B_r)\big)$$
is generated in degree $4$.
\end{thm}

I want to stress here that the constant $m$ can be calculated \emph{explicitly}, by tracing back through the results of \cite{CZ14,CL14} and calculating all the constants there precisely. The constant is large, but after one calculates it precisely, all the remaining steps in the algorithm presented in Sections \ref{sec:conedecomp} and \ref{sec:minmod} (apart from calculating the images of birational morphisms and finding suitable ample divisors) are then a mixture of convex geometry and linear algebra, which should be doable in modern computer algebra systems.

It is not completely obvious how to pass from Theorem \ref{t_efgs} to the Main Problem, so I explain this now. Let $(X,\Delta)$ be a projective $\Q$-factorial klt surface pair of log general type; passing to a log resolution, we may assume that $X$ is smooth. Pick some ample divisors $A_1',\dots,A_r'$ on $X$ and pick an element $A_0'$ of some non-empty linear system $|s (K_X+\Delta)|$ for $s>0$, which is possible since the divisor $K_X+\Delta$ is big. Since the cone of big divisors in $N^1(X)_\R$ is open, we can find distinct prime divisors $S_1,\dots,S_p$ on $X$ such that the numerical classes of all the divisors $A_0',\dots,A_r'$ are contained in the interior of the cone spanned by the numerical classes of $S_1,\dots,S_p$ in $\Div_\R(X)$. By replacing $A_i'$ by a divisor in the numerical class of $A_i'$ for each $i$, we may assume that $\Supp A_i'=\sum_{j=1}^pS_j$ for all $i$. Then take a large positive integer $q$ such that $\lfloor \frac1q A_i'\rfloor=0$ for all $i$. Set $A_i:=\frac1q A_i'$ for each $i$, and note that we have $K_X+\Delta+A_0\sim(\frac sq+1)(K_X+\Delta)$.

If the divisor $\Delta+\sum_{i=1}^p S_i$ has simple normal crossings support, then we can use Theorem \ref{t_efgs} to solve the Main Problem: namely, there exists a computable constant $m$ such that the ring
$$R\big(X;m(K_X+\Delta+A_0),m(K_X+\Delta+A_1),\dots,m(K_X+\Delta+A_r)\big)$$
is generated in degree $4$. Otherwise one has to pass to a log resolution of the pair $(X,\Delta+\sum_{i=1}^p S_i)$ and work there.

\begin{rem}
As pointed out to me by Frank-Olaf Schreyer, when $X$ is a \emph{smooth} surface of general type, a more efficient way to construct the minimal model of $X$ might be to apply \cite[Main Theorem]{Bom73} to first construct the canonical model of $X$ as the image of the pluricanonical map associated to the linear system $|5K_X|$, and then construct the minimal model of $X$ as the minimal resolution of the canonical model of $X$.
\end{rem}

\subsection*{The Main Problem in higher dimensions}

Inspired by \cite[Proposition 2.11]{CZ14} and \cite[proposition 2.18]{CL14}, I expect an analogue of Theorem \ref{t_efgs} to hold also on varieties $X$ in arbitrary dimension, by replacing $4$ by $\dim X+2$; the methods involved in the proof of \cite[Theorem 1.1]{CZ14} are probably crucial for the solution of the Main Problem on threefolds.

\section{Step 2: Cone decompositions}\label{sec:conedecomp}

The goal of the remainder of the paper is to demonstrate how once one solves the Main Problem, then one can \emph{in principle} be able to implement the MMP in a computer algebra system.

In order to present the second step, I will discuss briefly \emph{asymptotic geometric valuations}. This technical digression is not strictly necessary to understand the remainder of the algorithm, but it helps to understand the underlying philosophy. The details are in \cite{CL13} and in the survey \cite{Laz13}.

The finite generation of a divisorial ring $\mathfrak R$ has important consequences on the convex geometry of the cone $\Supp\mathfrak R$. We will relate the ring $\mathfrak R$ to the behaviour of linear systems $|D|$ for integral divisors $D\in\Supp\mathfrak R$ via asymptotic geometric valuations. 

Let $X$ be a smooth projective variety. A {\em geometric valuation\/} $\Gamma$ on $X$ is any valuation on $k(X)$ which is given by the order of vanishing at the generic point of a prime divisor on some birational model $Y\to X$, and we denote the value of this valuation on a $\Q$-divisor $D$ on $X$ by $\mult_\Gamma D$. If $D$ is moreover effective, then the {\em asymptotic order of vanishing\/} of $D$ along $\Gamma$ is 
$$o_\Gamma (D)=\inf\{\mult_\Gamma D'\mid D\sim_\Q D'\geq0\}.$$
Equivalently, if $\mult_\Gamma|kD|$ is the valuation at $\Gamma$ of a general element of the linear system $|kD|$, then 
$$o_\Gamma(D)=\inf\frac1k\mult_\Gamma|kD|$$
over all $k$ sufficiently divisible. 

The following result is \cite[Theorem 3]{CL13}: it is essentially contained in the proof of \cite[Theorem 4.1]{ELMNP}, and is in fact a mixture of commutative algebra and convex geometry. It provides the relation between the finite generation and the behaviour of linear systems, and yields an important decomposition of the cone $\Supp \mathfrak R$ into finitely many special rational polyhedral cones.

\begin{thm}
\label{thm:ELMNP}
Let $X$ be a smooth projective variety and let $D_0,\dots,D_r$ be $\Q$-divisors on $X$. Assume that the ring $\mathfrak R=R(X;D_0,\dots,D_r)$ is finitely generated. Then $\Supp \mathfrak R$ is a rational polyhedral cone and:
\begin{enumerate}[\normalfont (i)]
\item there is a finite rational polyhedral subdivision $\Supp \mathfrak R=\bigcup \mcal{C}_i$ into cones of maximal dimension, such that $o_\Gamma$ is linear on $\mcal{C}_i$ for every geometric valuation $\Gamma$ over $X$,
\item there exists a positive integer $k$ such that $o_\Gamma(kD)=\mult_\Gamma|kD|$ for every integral divisor $D\in \Supp \mathfrak R$.
\end{enumerate}
\end{thm}

As we will see in Section \ref{sec:minmod}, as soon one has the cone decomposition from Theorem \ref{thm:ELMNP}, one can calculate easily the outputs of the Minimal Model Program for a given projective klt pair of lof general type.

\medskip

So the main problem at this step of the algorithm is to construct such a decomposition. I sketch briefly the construction from \cite{ELMNP,ELMNPerratum}.

With notation from Theorem \ref{thm:ELMNP}, for each $m=(m_0,\dots,m_r)\in\N^{r+1}$ we denote by $\mathfrak{b}_m$ the ideal defining the base locus of the linear system $|m_0D_0+\ldots+m_rD_r|$. We obtain the \emph{$\N^{r+1}$-graded sequence of ideals} $\mathfrak{b}_\bullet := \{ \mathfrak{b}_m \}_{m \in \N^{r+1}}$ of ideal sheaves on $X$; that is, $\mathfrak{b}_0 = \OO_X$ and
\[  \mathfrak{b}_m \cdot \mathfrak{b}_{m'}  \subseteq \mathfrak{b}_{m + m'}\]
for all $m , m' \in \N^{r+1}$. Then we have the associated \emph{Rees algebra}
\[ R(\mathfrak{b}_\bullet) := \bigoplus_{m \in \N^{r+1}} \mathfrak{b}_m, \] 
which is finitely generated since the ring $\mathfrak{R}$ is finitely generated.

Then the arguments from \cite{ELMNP,ELMNPerratum} show that Theorem \ref{thm:ELMNP} follows from the following crucial proposition. 

\begin{pro}\label{pro:ELMNP} 
With notation as above, let $\mathcal C\subseteq \R^{r+1}$ be the convex cone spanned by $\N^{r+1}\subseteq\R^{r+1}$. Then there exist a smooth fan $\Delta$ with support $\mathcal C$ and a positive integer $d$ such that for every cone $\sigma\in\Delta$, if we denote by $e_1,\dots,e_s$ the generators of $\sigma\cap \N^{r+1}$,  then 
$$ \overline{\mathfrak{b}_{d{\sum_i}p_ie_i}}=\overline{\prod_i \mathfrak{b}_{de_i}^{p_i}},$$
for every $(p_1,\dots,p_s)\in\N^s$. Here $\overline{\mathfrak{b}}$ denotes the integral closure of an ideal $\mathfrak{b}$.
\end{pro}

Thus, the algorithm takes as an input the \emph{degrees of the generators of the ring $\mathfrak{R}$} and finds a \emph{rational polyhedral subdivision} of $\Supp\mathfrak{R}$ as an output.

The crucial point is this: the proof from \cite{ELMNPerratum} of Proposition \ref{pro:ELMNP} is convex geometric, thus finding a rational polyhedral subdivision of $\Supp\mathfrak{R}$ is implementable in suitable convex geometry software such as \texttt{Polymake} \cite{Polymake}. For applications to the MMP in Section \ref{sec:minmod} we actually need only Theorem \ref{thm:ELMNP}(i), and for this we do not need the full statement of Proposition \ref{pro:ELMNP}: following the proof in \cite[Proposition 1.1]{ELMNPerratum}, one first shows that each function $o_\Gamma$ can be written as the pointwise infimum of a family of linear functions -- these linear functions depend only on the multidegrees of some set of generators of $\mathfrak R$, and this is \emph{precisely the place} where the Main Problem is relevant. Once this is shown, one can find a subdivision as in Theorem \ref{thm:ELMNP}(i) by an explicit linear algebra argument, see \cite[Lemma 2.1 and Remark 2.3]{ELMNPerratum}.

For the sake of completeness, I mention here that Theorem \ref{thm:ELMNP}(ii) is also algorithmically implementable. Additional to Theorem \ref{thm:ELMNP}(i), there is only one crucial algebraic input, which corresponds to the case when $r=0$:

\begin{lem}\label{lem:Bourbaki}
If $R_0$ is a Noetherian ring and if $R=\bigoplus_{m\in\N}R_m$ is a finitely generated $R_0$-algebra, then there exists a positive integer $d$ such that $R_{dm}= R_d^m$ for every $m\in\N\setminus\{0\}$.
\end{lem}

The proof of this lemma in \cite[Chapter III, Section 1]{Bou98} is constructive: if one knows the degrees of some generators of $R$, then one can calculate $d$ explicitly. 

Therefore, one should be able to implement Theorem \ref{thm:ELMNP} in a computer algebra system (in the sense mentioned earlier: take as input the degrees of the generators of the ring $\mathfrak{R}$ and get as output a rational polyhedral subdivision of $\Supp\mathfrak{R}$), although improvements and simplifications to the algorithm from \cite{ELMNPerratum} have to be made, since the calculations can get expensive quickly as the number of generators of $\mathfrak R$ or their multidegrees grow. 

Once Theorem \ref{thm:ELMNP} is implemented, we apply it for a projective $\Q$-factorial klt pair $(X,\Delta)$ of log general type. We let $A_1,\dots,A_r$ be ample divisors on $X$, and set $D_0:=K_X+\Delta$ and $D_i:=K_X+\Delta+A_i$ for $i=1,\dots,r$. Then, as explained above, as soon as we have a solution to the Main Problem for the ring $\mathfrak R=R(X;D_0,\dots,D_r)$ (for instance, either by finding generators of $\mathfrak R$ explicitly in explicit examples, or by knowing a general upper bound for their multidegrees by general theoretical arguments as in \cite{CZ14,CL14}), we can calculate the cone subdivision as in Theorem \ref{thm:ELMNP}(i) in a computer algebra system. As mentioned in Section \ref{sec:graded}, the biggest progress on this problem is currently on surfaces, where one can calculate a bound explicitly by Theorem \ref{t_efgs}.

\section{Step 3: Minimal models}\label{sec:minmod}

In the final step, we will see that the cone decomposition from the previous section applied to a particular choice of a divisorial ring on a klt pair $(X,\Delta)$ of log general type immediately provides the steps in a Minimal Model Program of $(X,\Delta)$, as well as a minimal model of $(X,\Delta)$.

So let $(X,\Delta)$ be a projective $\Q$-factorial klt pair of log general type and let $A_1,\dots,A_r$ be ample divisors on $X$ whose classes generate $N^1(X)_\R$. Thus, we may assume that $r$ is the Picard rank of $X$. Moreover, by the Cone theorem, possibly by replacing $A_i$ by $(2\dim X+1)A_i$, we may assume from the start that each divisor $D_i:=K_X+\Delta+A_i$ is ample for $i=1,\dots,r$: this will make finding a divisor $H$ below easier, but is not essential.

Then by Theorem \ref{thmA} the ring $\mathfrak R=R(X;K_X+\Delta,D_1,\dots,D_r)$ is finitely generated. Therefore, by Theorem \ref{thm:ELMNP} the cone $\mathcal C:=\Supp \mathfrak R$ has a finite rational polyhedral subdivision 
$$\mathcal C=\bigcup_{i=1}^p \mcal{C}_i$$
into cones of maximal dimension, such that the asymptotic valuation function $o_\Gamma$ is linear on each $\mcal{C}_i$ for every geometric valuation $\Gamma$ over $X$. By subdividing $\mcal C$ further, we may assume that each cone $\mathcal C_i$ is contained in one of the two half-spaces of the vector space $\R (K_X+\Delta)+\sum\R D_i$ bounded by each hyperplane which contains a face of some $\mathcal C_j$, where $1\leq i,j\leq p$. Note that the cone $\mcal C$ is equal to the cone $\R_+(K_X+\Delta)+\sum\R_+D_i$, since $K_X+\Delta$ is big.

Pick an ample divisor $H$ on $X$ such that, denoting 
$$I:=\big\{t(K_X+\Delta)+(1-t)H\mid 0\leq t\leq1\big\},$$
then we have the following property: if $I$ intersects some cone $\mathcal C_i$, then it intersects its interior -- it is clear that one can easily achieve this. By relabelling, we may assume that $I$ intersects precisely the cones $\mathcal C_1,\dots,\mathcal C_k$ for some $k\leq p$, so that
$$K_X+\Delta\in\mathcal C_k, \quad H\in\mathcal C_1\quad\text{and}\quad\mathcal C_i\cap\mathcal C_{i+1}\cap I\neq\emptyset\text{ for all }i=1,\dots,k-1.$$

Then we have:

\begin{thm}\label{thm:minimal}
With notation as above, let $D$ be any divisor in the interior of $\mathcal C_k$. Then the variety $\Proj R(X,D)$ is a minimal model of $X$.
\end{thm}

This follows from the proof of \cite[Theorem 6]{CL13} and from the proof of a more general statement \cite[Theorem 5.4]{KKL16}. I briefly sketch the main steps of the proof, which will in particular explain how to run the MMP with scaling of $H$ in this context. This is important when one is interested in all the steps of the MMP instead of just in its main result.

The proof is by induction on $k$. Denote by $\pi\colon \Div_\R(X)\to N^1(X)_\R$ the natural projection. Then one can show that the cone $\mcal C\cap \pi^{-1}\bigl(\Nef (X)\bigr)$ is the union of some of the cones $\mathcal C_i$, where $1\leq i\leq p$. In particular, we may assume that there exists $1\leq k_0\leq k$ such that 
$$\mathcal C_i\subseteq \pi^{-1}\bigl(\Nef (X)\bigr)\quad\text{for }i=1,\dots,k_0$$
and
$$\mathcal C_i\not\subseteq \pi^{-1}\bigl(\Nef (X)\bigr)\quad\text{for }i=k_0+1,\dots,k.$$
Pick any divisor $G$ in the interior of $\mathcal C_{k_0+1}$ and set $X':=\Proj R(X,G)$. Then the proof of \cite[Theorem 5.2]{KKL16} shows that the induced birational map $f\colon X\dashrightarrow X'$ is precisely the first step in the $(K_X+\Delta)$-MMP with scaling of $H$.

Define rational polyhedral cones $\mathcal C_i':=f_*\mathcal C_i\subseteq\Div_\R(X')$ and denote by $\pi'\colon \Div_\R(X')\to N^1(X')_\R$ the natural projection. Then one can again show that there exists $k_0+1\leq k_1\leq k$ such that 
$$\mathcal C_i'\subseteq (\pi')^{-1}\bigl(\Nef (X')\bigr)\quad\text{for }i=k_0+1,\dots,k_1$$
and
$$\mathcal C_i'\not\subseteq (\pi')^{-1}\bigl(\Nef (X')\bigr)\quad\text{for }i=k_1+1,\dots,k.$$
Pick any divisor $G'$ in the interior of $\mathcal C_{k_1+1}$ and set $X'':=\Proj R(X,G')$. Then the induced birational map $f'\colon X'\dashrightarrow X''$ is precisely the second step in the $(K_X+\Delta)$-MMP with scaling of $H$. By continuing this procedure we realise the whole MMP with scaling.

The procedure just described calculates the maps $f,f',\dots$ in the MMP above, but there is some work involved in finding the indices $k_j$. Moreover, one calculates each subsequent map on a different variety: for instance, $f'$ was determined by data on $X'$ and not on $X$ itself. One can, however, proceed also as follows: For each $i=1,\dots,k$ one can pick a divisor $G_i$ in the interior of $\mathcal C_i$ and set $X_i:=\Proj R(X,G_i)$. Then one obtains birational maps
$$X=X_1\dashrightarrow X_2\dashrightarrow\ldots\dashrightarrow X_k,$$
where $X_k$ is a minimal model of $X$. Some of these birational maps will actually be isomorphisms, hence there will be some redundancy in the process. However, the clear advantage is that one has to calculate the varieties $X_i$ on the variety $X$ itself, which saves computer power.

\section{Conclusion}

An algorithm for an implementation of the Minimal Model Program in a computer algebra system was proposed in Sections \ref{sec:graded}, \ref{sec:conedecomp} and \ref{sec:minmod}, for pairs $(X,\Delta)$ of log general type.

The algorithm proceeds in three steps. In Step 1 one considers a carefully chosen multigraded ring $\mathfrak R$ on $X$. The \emph{Main Problem} in the algorithm is to bound the multidegrees of its generators. This problem has a solution on surfaces, therefore the MMP on surfaces should be implementable on a computer in the near future. Bounding the multidegrees of generators of such graded rings $\mathfrak R$ in higher dimensions is a topic of future theoretical research.

Once the Main Problem is implemented, the remainder of the algorithm becomes a mixture of linear algebra and convex geometry. In Step 2 one uses as an \emph{input} the information on the multidegrees of generators of $\mathfrak R$ as in the Main Problem to find a suitable decomposition of the support of $\mathfrak R$ into finitely many rational polyhedral cones with special properties. In Step 3, from this decomposition one can immediately obtain a minimal model of $(X,\Delta)$ and the steps in a $(K_X+\Delta)$-MMP with scaling of a certain ample divisor. All three steps can be implemented independently of each other.
	
	\bibliographystyle{amsalpha}
	\bibliography{biblio}
	
\end{document}